\documentclass[11pt]{article}
\usepackage{graphicx}
\usepackage{amsmath}
\usepackage{mathrsfs,amsfonts,bbm}
\usepackage{amsfonts}
\usepackage{amsmath,amssymb,amsfonts}
\usepackage{amssymb}
\usepackage{color}
\definecolor{c}{rgb}{0.9,0.3,0.1}
\definecolor{b}{rgb}{0.1,0.3,0.9}

\usepackage{ulem}

 \allowdisplaybreaks

\newtheorem{theorem}{Theorem}[section]

\newtheorem{lemma}[theorem]{Lemma}

\newtheorem{definition}[theorem]{Definition}

\newtheorem{remark}[theorem]{Remark}

\newtheorem{hypothesis}[theorem]{Hypothesis}

\numberwithin{equation}{section}

\def\al{{\alpha}}\def\de{{\delta}}

\def\si{{\sigma}}
\def\<{\left<}\def\>{\right>}\def\({\left(}\def\){\right)}

\newfam\msbmfam\font\tenmsbm=msbm10\textfont
\msbmfam=\tenmsbm\font\sevenmsbm=msbm7
\scriptfont\msbmfam=\sevenmsbm\def\bb#1{{\fam\msbmfam #1}}

\def\RR{\bb R}
\newcommand{\E}{{\mathbbm{E}}}

\def\cF{{\cal F}}
\def\cM{{\cal M}}

\newcommand{\EE}{{\mathbb{E}}}
\newcommand{\PP}{{\mathbb{P}}} 

\newcommand{\R}{{\mathbf{R}}}

\newcommand{\bP}{{\mathbbm{P}}}

 \def\qed{\hfill $\Box$}

\def\sF{\mathcal F}

\def\wt{\widetilde}

\def\eps{\varepsilon}

\def\R{\mathbbm{R}}
\def\1{\mathbbm{1}}

\def\1{\mathbbm{1}}

\def\pf{\noindent {\bf Proof.} }

\numberwithin{equation}{section}
\begin{document}
\title{Fully coupled forward-backward stochastic differential equations driven by sub-diffusions
 }
  \author{ Shuaiqi Zhang\footnote{Research partially  supported  by the Fundamental Research Funds for the Central Universities (Grant No. 2020ZDPYMS01) 
 and Nature Science Foundation of Jiangsu Province (Grant No. BK20221543)  and National Natural Science Foundation of China (Grant No. 12171086).  }   \qquad \hbox{and} \qquad 
 {Zhen-Qing Chen\footnote{ Corresponding author. Research  partially supported by a Simons Foundation fund. } 
 }    }
\maketitle
\date
\begin{abstract}   
	In this paper,     we   establish the existence and uniqueness of  fully coupled  forward-backward stochastic differential equations
	(FBSDEs in short)   driven by  anomalous sub-diffusions $B_{L_t}$ under suitable monotonicity conditions on the coefficients.  
Here $B$ is a Brownian motion on $\R$ and $L_t:= \inf\{r>0: S_r>t\}$, $t\geq 0,$
   is  the inverse of  a subordinator  $S$ with drift $\kappa >0$ that is independent of $B$.   
Various   a   priori estimates on the solutions of the FBSDEs are also presented.   
\end{abstract}

 \medskip
 

\medskip


\section{Introduction}\label{S1}

 Forward and backward stochastic differential equations (FBSDEs) driven by Brownian motion
 is an important and effective tool in stochastic control theory and its applications including mathematical finance. They arise naturally as the adjoint equation for the stochastic Pontryagin-type maximum principle \cite{Peng1993}, and in hedging options model for large investors \cite{CM}.
 They played an instrumental role in the resolution of  the Black's consol rate conjecture \cite{DMP}.

Sub-diffusions  are a class of random  processes that  describe the motion of particles that move slower than Brownian motion,
for example, due to particle sticking and trapping. 
It has been observed that many phenomena in real world systems such as porous media, biological systems and financial markets can be and should be modeled by sub-diffusion.  For example, it can be used to model the spread of diseases in populations, and transport of particles through
porous media such as soils.
  Prototypes of anomalous sub-diffusions are Brownian motions time-changed by the inverse of subordinators; see, e.g., 
\cite{Meerschaert2004, MK}.    They do not have Markov properties.

The goal of this paper is to study fully coupled forward-backward stochastic differential equations  (FBSDEs) driven by non-Markov anomalous sub-diffusions. We are led to this study from several  aspects. 
Stochastic maximal principles (SMPs in short) 
for  backward stochastic differential equations (BSDEs in short) driven by sub-diffusion have recently been investigated in \cite{ZhangChen2023}, see also \cite{NN} for a sufficient SMP.   As pointed out in \cite{ZhangChen2023}, 
among other things, SDEs  driven by sub-diffusions can be used to model financial activities in bear markets 
when  activities are much slower. 
It is reasonable to expect that the adjoint equation for the   SMP  for BSDEs 
recently obtained in \cite{ZhangChen2023}
can be described by a suitable FBSDE driven by a sub-diffusion. The  hedging options model for large investors 
 and the  consol rate problem in bear markets  
  also give a motivation to study FBSDEs driven by sub-diffusions.

Existence and uniqueness of solutions  to  coupled FBSDEs driven by Brownian motion has been extensively studied in literature.
Motivated by the  stochastic differential utility in \cite{DE} and the seminal work of Pardoux and Peng \cite{PP} on BSDEs, 
 Antonelli   \cite[Theorem 3.1]{Antonelli1993}  studied the existence and uniqueness of strongly coupled FBSDEs driven by Brownian motion
 and showed that,  under  some smallness Lipchitz conditions, 
 they have unique solution  in  small time duration   using a contraction map argument.
  Furthermore, a counterexample   \cite[Exampe 1]{Antonelli1993}  is given 
   showing that Lipschitz condition is not enough for the existence of solution to FBSDE over general bounded time intervals.   
   In    Ma-Protter-Yong \cite{MPY},  a four step scheme is developed that gives existence and uniqueness of coupled FBSDEs over
   arbitrary time duration via a quasi-linear partial differential equation but requiring the non-degeneracy of the forward SDE
   and non-randomness of the coefficients. 
    In  Hu and Peng \cite{HuPeng1995}, existence and uniqueness results of fully coupled FBSDEs driven by Brownian motion
over any   prescribed  time duration are obtained under some monotonicity conditions using a method of continuation.
This result has been extended further in Yong \cite{Y97}, Peng and Wu \cite{PengWu1999} and Hu \cite{Hu2000}.
We refer the reader to the two excellent books \cite{Ma1999, Zhang2017} and the references therein for more information on the theory of 
FBSDEs driven by    Brownian    motion. For further applications of FBSDEs, see \cite{CCM}.

To the best knowledge of the authors, fully coupled forward-backward stochastic differential equations  (FBSDEs) driven by non-Markov anomalous sub-diffusions have not been previously studied. In this paper, we adapt the method of continuation to establish the existence and uniqueness 
for fully coupled FBSDEs with random coefficients over arbitrary time duration. 
We point out that allowing the coefficients to be random is crucial for stochastic control of FBSDEs.
We further derive a priori estimates on the solutions of FBSDEs, which will play an important role in our study of stochastic optimal control 
of FBSDEs in \cite{ZC2}. Subdiffusions can stand still over infinitely many small time intervals. On the other hand, 
sub-diffusions have Brownian motion as their extreme case.
 The main difference between  FBSDEs  driven by Brownian motion and by sub-diffusions is that the latter ones 
  can be degenerate in the sense that sub-diffusions can remain constant during infinitely many random time intervals.

\medskip

 Suppose that $S=\{S_t; t\geq 0\}$ is a subordinator with drift $\kappa \geq 0$ and  
L\'evy measure $\nu$; that is, $S_t=\kappa t+ S_t^0$, where $S_t^0$ is a driftless subordinator with
  L\'evy measure $\nu$.
 Let $L$ be the inverse of $S$, that is, 
 $$
 L_t = \inf\{r>0: S_r >t\}, \quad t\geq 0.
 $$
 The inverse subordinator $L_t$ is continuous in $t$  but stays constant during infinitely many time periods which are resulted from the infinitely many jumps by the subordinator $S_t$ when its L\'evy measure is non-trivial.  
Let $B$ be a  Brownian motion independent of $S_t$. 
  Its time-change process $B_{L_t}$ by the inverse subordinator $L_t$ is called a sub-diffusion $B_{L_t}$,
which  is a continuous martingale with
$\langle B_{L_t} \rangle = L_t$  but is not a Markov process. 
   Note that  $B_{L_t}$ stays flat during the time periods when $L_t$ stays constant. 
   For any $\eps >0$,  the jumps of $S_t$ of size larger than $\eps$ occurs 
   according to a Poisson process with parameter $\nu (\eps, \infty)$.  When the L\'evy measure $\nu$ of 
   the subordinator $S$ is infinite, 
   during any fixed time intervals, 
   $L_t$ has infinitely many small time periods  
      but  only finite many time intervals larger than $\eps$  during which it stays constants
   Thus the sub-diffusion $B_{L_t}$  matches well with the phenomena that the financial market constantly  has small corrections 
    but long bear market occurs only sporadically.

 \medskip 

 In this paper, we consider  fully coupled FBSDEs driven anomalous sub-diffusions  of the following form
  \begin{equation*} 
\left\{\begin{aligned}
dx (t)=&\ b(\omega, t, x (t), y(t)  )dt+\delta (\omega, t, x(t),  y(t),z(t) )d L_{(t-a)^+ } + \sigma   (\omega, t, x(t),  y(t) ,z(t) )  dB_{L_{(t-a)^+ }},\\
 - dy (t)=&\  g(\omega, t, x (t), y(t) )dt  +  h (\omega, t, x(t), y(t), z(t) )d L_{(t -a)^+} -  z (t)dB_ {L_{(t -a)^+}}  ,\\
x (0)=&\ x_{0}  \quad \hbox{and} \quad y (T)=\phi(\omega, x (T) ).
\end{aligned}
\right.
\end{equation*}
  Here $0\leq a<T<\infty$, where $a\geq 0$ models the initial wake up time for  the sub-diffusion
$B_{L_{(t-a)^+}}$ to become active. 
See \eqref{FBSDE} below for details. 
 Suppose that the subordinator $S$ has positive drift $\kappa$.
We establish the existence and uniqueness 
for fully coupled FBSDEs with random coefficients over arbitrary time duration using a method of continuation 
and the martingale representation theory for sub-diffusions recently established in \cite{ZhangChen2023}. 
  We further  derive various a priori estimates on the solutions of FBSDEs that 
will play a crucial role in the study \cite{ZC2} of their control problems.
 An interesting feature of FBSDE studied in this paper is that there are two time  scales $ dt$ and $dL_t$.  
 The   FBSDEs driven by sub-diffusions exhibit  a combined deterministic and stochastic features.  
  Note that when the L\'evy measure $\nu$ for the subordinator  vanishes,  $S_t=\kappa t$ and so 
$B_{L_t}= B_{t/\kappa}$ reduces to a Brownian motion. Thus 
the results in this paper  not only  recover but also extend in a ``continuous way" the corresponding results in
the classical Brownian setting.

 \medskip
 
   For notational convenience, all the processes in this paper are assumed to be real-valued. However,
   the results in this paper work for multi-dimensional case with straightforward modification. 
   Throughout  this paper, we use notation $:=$ as a way of definition. 
  For a stochastic process $x$, we use the notation $x(t)$ and $x_t$ interchangeably, 
  to denote its state or position  at time $t$.

 \medskip
 
 The rest of this paper is organized as follows. In Section \ref{S2}, we recall some facts from \cite{ZhangChen2023} about  inver subordinators and 
 anomalous sub-diffusions that will be used in this paper. 
 Existence and uniqueness of the solution to fully coupled FBSDEs with random coefficients 
  over any time duration  is given   in Theorem \ref{exsi-uniqu} of Section \ref{S3} 
  under a monotonicity condition by a method of continuation.  
  Moreover,   we derive a priori estimates for the unique solution $(x(t), y (t), z (t))$ of the FBSDE.

\section{Preliminary about Sub-diffusion}\label{S2}

In this section, we recall some results from  \cite{ZhangChen2023} that will be used later in this paper.
Although the sub-diffusion itself is not a Markov process, we can make it Markov by adding an auxiliary overshoot process. 
 
 \begin{theorem}\label{T:3.1} 
  Suppose that $B$  is a  standard Brownian motion on $\RR$ starting from $0$,  $S$ is any subordinator that is independent of $B$ with $S_0=0$,
and $L_t:=\inf\{r>0:S_r>t\}$.
Then
\begin{equation} \label{e:2.1}
\wt X_t:=(X_t, \, R_t):= \left(x_0+ B_{L_{(t-R_0)^+}},  \, R_0+ S_{L_{(t-R_0)^+}} - t  \right), \quad t\geq 0,
\end{equation}
with $\wt X_0= (x_0, R_0) \in \RR \times [0, \infty)$   is a time-homegenous Markov process taking values in $\RR\times [0, \infty)$.
\end{theorem}

\bigskip

Note that for each fixed $t>0$,  $S_{L_t}>t$ happens with positive probability.
 On $\{S_{L_t}>t\}$, the inverse local time $L_s$ and, consequently, 
 the sub-diffusion $B_{L_s}$ remain flat during the time interval $[t, S_{L_t} ]$. 
 We call $R_t:=R_0+ S_{L_{(t-R_0)^+}} - t $ an overshoot process with initial value $R_0$. 
 It measures how much time it would take for the anomalous sub-diffusion $X_t:= x_0+ B_{L_{(t-R_0)^+}}$ to wake up   from time $t$. 
 
\medskip

  Let  $\{\sF'_t\}_{t\geq 0}$ be the minimum augmented filtration generated by $ X$
  and    $\{\sF_t\}_{t\geq 0}$ be     the minimum augmented filtration generated by 
  $\wt X=(X, R)$.   Clearly, $\sF'_t\subset \sF_t$ for every $t\geq 0$. 
Note that process $\wt X_t=(X_t, R_t)$  
  depends on the initial $a\geq 0 $ of $R_0$, so do the filtrations $\{\sF_t\}_{t\geq 0}$ and $\{\sF'_t\}_{t\geq 0}$. When $R_0=a$ for some deterministic $a\geq 0$, for emphasis, sometimes we denote them
by $\wt X^a_t=(X^a_t, R^a_t)$, $\{\sF^a_t\}_{t\geq 0}$ and $\{ {\sF^{a}_t}'\}_{t\geq 0}$, respectively.  Clearly, $\sF^a_t$ and $ {\sF^{a}_t}' $ are trivial for $t\in [0, a]$.

\medskip

In the rest of this paper, we assume the subordinator $S$ has positive drift $\kappa >0$. 
In this case,  for any $t, s>0$,
 $$
 0\leq L_{t+s} -L_t \leq s/\kappa.
 $$
 So almost surely
\begin{equation}\label{e:2.2}
 \frac{dL_t}{dt}    \ \hbox{ exists} \quad \hbox{with} \quad 0\leq  \frac{dL_t}{dt}  \leq 1/\kappa   \    \hbox{ for a.e. }   t>0.
\end{equation} 

\section{Existence and uniqueness for FBSDEs driven by sub-diffusion}\label{S3}

Fix  $0\leq a <T<\infty$. 
Consider the following coupled 
FBSDE driven by sub-diffusion $B_{L_t}$ for $t\in [0, T]$:
  \begin{equation}\label{FBSDE}
\left\{\begin{aligned}
dx (t)=&\ b(\omega, t, x (t), y(t)  )dt+\delta (\omega, t, x(t),  y(t),z(t) )d L_{(t-a)^+ } + \sigma   (\omega, t, x(t),  y(t) ,z(t) )  dB_{L_{(t-a)^+ }},\\
 - dy (t)=&\  g(\omega, t, x (t), y(t) )dt  +  h (\omega, t, x(t), y(t), z(t) )d L_{(t -a)^+} -  z (t)dB_ {L_{(t -a)^+}}  ,\\
x (0)=&\ x_{0}  \quad \hbox{and} \quad y (T)=\phi(\omega, x (T) ), 
\end{aligned}
\right.
\end{equation}
where for each fixed $x, y, z\in \R$, $b(\omega, t, x, y)$ and $g(\omega, t, x, y)$ are real-valued $\{\sF'_t\}$-progressively measurable random processes defined on $\Omega \times  [0, T]$,
$\delta (\omega, t, x, y, z)$,  $\si(\omega, t, x, y, z)$ and  $h (\omega, t, x, y, z)$ are $\{\sF'_t\}$-progressively measurable random processes defined on  $\Omega \times  [0, T]$, 
 $ \phi (\omega, x) $ is an $L^2$-integrable $\sF_T'$-measurable real-valued random variable defined on $\Omega$.
 For notational simplicity, we will typically drop $\omega$ from the expressions  of the above random processes or variables.

 \medskip

 \begin{hypothesis}\label{HP0}
 \begin{enumerate}
 \item[\rm (i)] 
  $ \delta,  \, \si$ and  $h $ are uniformly Lipschitz continuous in $(x,y,z)$,    
  $b$ and $g$ are uniformly Lipschitz continuous in $(x,y)$, 
  and $\phi$ is uniformly Lipschitz continuous in $x$ with Lipschitz constant $L\geq 1$
   with
  $$
  \EE \left[ \delta^2 ({\bf 0}) +\si^2  ({\bf 0})+h^2 ({\bf 0})   + b^2 ({\bf 0})+g^2 ({\bf 0}) + \phi^2(0)\right] <\infty.  
	$$
	Here ${\bf 0}$ denotes the origin in $\R^4$ or in  $\R^3$.

 \item[\rm (ii)] 
There is some $c>0$  so that the following two monotonicity conditions hold
$\PP$-a.s.: 
for any $t>0$,  $x_1, x_2, y_1, y_2, z_1$ and $z_2$ in $\R$: 
\begin{eqnarray}\label{m1}
  &&  \big(  b(t,x_1,y_1 )  -b(t,x_2,y_2  ) \big) (y_1-y_2  ) 
    -   \big(  g(t,x_1,y_1 )  -g(t,x_2,y_2 ) \big) (x_1-x_2  ) \nonumber\\
 &&  \quad \leq -c \big( | x_1-x_2 |^2  + | y_1-y_2 |^2   \big) ,
\end{eqnarray}
and
\begin{eqnarray}\label{m2}
&&    \big(  \si(t,x_1,y_1,z_1)  -\si(t,x_2,y_2,z_2) \big) (z_1-z_2  ) +  \big(  \de(t,x_1,y_1,z_1)  -\de(t,x_2,y_2,z_2) \big) (y_1-y_2  )   \nonumber\\
&& \qquad \  -\big(  h(t,x_1,y_1,z_1 )  -h(t,x_2,y_2,z_2 ) \big) (x_1-x_2  ) \nonumber\\ 
&&\quad \leq 
  -c \big( | x_1-x_2 |^2  + | y_1-y_2 |^2 + | z_1-z_2 |^2  \big).
\end{eqnarray}

\item[\rm (iii)] The function $\phi (x)$ is non-decreasing in $x$.
    \end{enumerate}
  \end{hypothesis}
  
\medskip

\begin{remark} \label{R:3.2} \rm  
\begin{enumerate}  
\item[\rm (i)]  Suppose \eqref{m1} holds. Taking $y_1=y_2$ in \eqref{m2} yields that  $\PP$-a.s., 
\begin{equation}\label{e:3.4}
\big(  g(t,x_1,y )  -g(t,x_2, y ) \big) (x_1-x_2  )   \geq c  |x_1-x_2|^2
\quad \hbox{for any } t>0, x_1, x_2, y\in \R.
\end{equation}  
Similarly, taking $x_1=x_2$ in \eqref{m1} yields that $\PP$-a.s., 
\begin{equation}\label{e:3.5} 
\big(  b(t,x,y_1 )  -b(t,x,y_2  ) \big) (y_1-y_2  ) \leq - c |y_1-y_2|^2
\quad \hbox{for any } t>0, x, y_1, y_2\in \R.
\end{equation}

\item[\rm (ii)]  Conversely, suppose that $g(t, x, y)$ and $b(t, x, y)$ are functions that satisfy conditions
\eqref{e:3.4} and \eqref{e:3.5}, respectively, and that $b$ is uniformly Lipschitz continuous in
$x$ with Lipschitz constant $L$ and $g$ is uniformly Lipschitz continuous in
$y$ with Lipschitz constant $L$. If $L<c/2$, then \eqref{m1} holds with $c/2$ in place of $c$.  
This is because 
\begin{eqnarray*}
&& (  b(t,x_1,y_1 )  -b(t,x_2,y_2) ) (y_1-y_2  )  -  (  g(t,x_1,y_1)  -g(t,x_2,y_2)) (x_1-x_2)  \\
&=& (  b(t,x_1,y_1)  -b(t,x_1,y_2) ) (y_1-y_2) 
+ (b(t,x_1,y_2) - b(t,x_2,y_2) )(y_1-y_2)  \\
&& - (  g(t,x_1,y_1)  -g(t,x_2,y_1 )) (x_1-x_2)
  -  (  g(t,x_2,y_1)  -g(t,x_2,y_2 )) (x_1-x_2  )  \\
&\leq & -c |y_1-y_2|^2 + L  |x_1-x_2|\, |y_1-y_2| - c |x_1-x_2|^2 + L |x_1-x_2|\, |y_1-y_2| \\
&\leq & -c (|x_1-x_2|^2 + |y_1-y_2|^2) + L (|x_1-x_2|^2 + |y_1-y_2|^2) \\
&\leq & -(c/2) (|x_1-x_2|^2 + |y_1-y_2|^2) .
\end{eqnarray*}
\end{enumerate}
Similar remark applies  condition \eqref{m2} as well.   \qed 
\end{remark}

\medskip

Suppose $0\leq a <T<\infty$. We  define a Banach norm 
$\| \cdot  \|_{ \mathcal M_a  [0,T]}$  on the space  
  \begin{eqnarray*}
\mathcal M_a [0,T] 
&:=& \Big\{\Theta(t)=  (x(t), y(t), z(t)): \hbox{ $\Theta (t)$ is an $\{\sF'_t\}_{t\in [0, T]}$-progressively measurable process}  \\
&& \hskip 0.3truein 
 \hbox{ on $[0, T]$ with } 
  \EE  \Big[  |x(0)|^2+  \int_0^T  (|x(t)|^2+  |y(t)|^2 ) dt  +  \int_0^T    | z (t)|^2   dL_ {(t-a)^+ }  \Big] <\infty \Big\} 
  \end{eqnarray*}
by 
\begin{eqnarray}
\|   \Theta  \|_{ \mathcal M_a  [0,T]}  
:=  
\left(    \EE  \left[ |x(0)|^2+  \int_0^T  (|x(t)|^2+  |y(t)|^2 ) dt   +  \int_0^T   | z (t)|^2    dL_ {(s-a)^+ }  \right]  \right)^{1/2}.
\end{eqnarray}

 \medskip
 
  Note that under the assumption that $S$ has positive drift $\kappa >0$, we have by \eqref{e:2.2} that 
 norm $\|   \Theta  \|_{ \mathcal M_a  [0,T]}  $ is comparable to
 $$
 \left(    \EE  \left[ |x(0)|^2+  \int_0^T  (|x(t)|^2+  |y(t)|^2 ) dt   +  \int_0^T   | \Theta  (t)|^2    dL_ {(s-a)^+ }  \right]  \right)^{1/2}.
 $$

\medskip 

\begin{definition}\label{XYZsolution}
Let $0\leq a <T<\infty$. 
A triple of processes $ \Theta (\cdot):= (x(\cdot), y(\cdot), z(\cdot)) \in \mathcal M_a [0,T]$ is called an adapted solution of (\ref{FBSDE}) if  for any $t\in [0, T]$,
   \begin{equation*}\label{FBSDE-def}
\left\{\begin{aligned}
x (t)=&x_0 + \int_0^t    b(s, \Gamma (s) )ds+  \int_0^t  \delta (s,\Theta (s) )d L_{(s-a)^+ } +  \int_0^t  \sigma   (s, \Theta (s) )  
dB_{L_{(s-a)^+ }},\\
  y (t)=&\phi(x (T) ) + \int_t^T   g(s, \Gamma (s) )ds  +\int_t^T  h (s, \Theta (s) )d L_{(s-a)^+}  -\int_t^T  z (s)dB_ {L_{(s -a)^+}}  ,
\end{aligned}
\right.
\end{equation*}
where $\Gamma (s):=(x (s), y(s)) $. 
We say the solution to \eqref{FBSDE} is unique if $ (x(\cdot), y(\cdot), z(\cdot)) \in \mathcal M_a[0,T]$ and
 $(\wt x(\cdot), \wt y(\cdot), \wt z(\cdot)) \in \mathcal M_a [0,T]$ are two solutions of \eqref{FBSDE},
  then $\wt x_t=x_t$ and   $\wt y_t=y_t$ for all $t\in [0, T]$ with probability one,
   and $\EE \int_0^T |z_s -\wt z_s|^2 dL_ {(s-a)^+ }    =0$.
  \end{definition}

\begin{remark} \rm 
By an argument similar to that of \cite[Lemma 4.4]{ZhangChen2023},
under  the Lipschitz assumption in  Hypothesis \ref{HP0}(i),  one can deduce that there is a constant $C>0$ so that 
for any adapted solution $\Theta (\cdot):= (x(\cdot), y(\cdot), z(\cdot)) $  of (\ref{FBSDE}), 
$$
  \EE  \left[  \sup_{t\in [0, T]} ( |x(t)|^2+   |y(t)|^2 )   +  \int_0^T   | z (t)|^2    dL_ {(s-a)^+ }  \right] 
  \leq C      \left( |x_0|^2 +   \EE  \left[     |\phi (0)|^2 \right] + \| \Theta    \|^2_{ \mathcal M_a  [0,T]} \right) .    
$$
  Moreover, by \eqref{e:2.2}, 
\begin{eqnarray}
&&   \EE  \left[  \sup_{t\in [0, T]} ( |x(t)|^2+   |y(t)|^2 )   +  \int_0^T   | \Theta  (t)|^2    dL_ {(s-a)^+ }  \right]   \nonumber \\
&  \leq&  C(T, \kappa) \EE  \left[  \sup_{t\in [0, T]} ( |x(t)|^2+   |y(t)|^2 )   +  \int_0^T   | z (t)|^2    dL_ {(s-a)^+ }  \right] .
\end{eqnarray} 
 
\qed 
\end{remark}

\medskip

 In this section,  we will adapt the method of continuation from  \cite{HuPeng1995, PengWu1999}
to establish the following existence and uniqueness of solutions to FBSDE \eqref{FBSDE}  with random coefficients  
in any prescribed time duration.  It is  the main result of this paper.

\medskip

\begin{theorem}\label{exsi-uniqu}
Under Hypothesis \ref{HP0},  for any $0\leq a<T<\infty$,  
the  FBSDE \eqref{FBSDE} admits a unique solution $\Theta (t):=(x(t), y(t), z(t))$ in $\mathcal M_a [0, T]$.
  Moreover, there is a constant $C_0= C_0(T,  \kappa, L, c)>0$, where $L$ and $c$ are the positive constants in Hypothesis \ref{HP0}, so that 
\begin{eqnarray} \label{e:3.7a}
  && \EE  \left[ \sup\limits_{[0,T]}  \( x(t)^2+y(t)^2  \)    +\int_0^T   |z (t) |^2  d L_{(t-a)^+ }  \right] \nonumber\\
  & \leq &   C_0 \bigg( x_0^2 +      \EE \left[          \phi (0)^2 \right] +  \EE  \int_0^T  \left(
   b (  t,0,0)^2     +  g  (  t,0,0)^2  \right) dt  \\
           && \quad + \EE \int_0^T \left(    \de (  t,0,0,0  )^2    + h ( t,0,0,0)^2    +
      \si ( t,0,0,0  )^2   \right)dL_{(t-a)^+}       \bigg) .   \nonumber 
          \end{eqnarray}
\end{theorem}
 
 We will establish the above theorem through two key  lemmas.

 \begin{lemma} \label{L:3.4} 
 Let $T>0$.
Suppose that  $\phi_0 \in L  ^2 (\sF'_T)$, and $b_0, \de_0,  \,\si_0 , \, g_0$ and  $h_0$ are $\{\sF'_t\}_{t\in [0, T]}$- progressively measurable processes on $\Omega \times [0, T]$ 
  with
$$
\EE \left[  \int_0^T  \( b_0(s)^2 +g_0 (s)^2 \) ds  + \int_0^T \( \de_0(s)^2 + h_0(s)^2 + \si_0^2   \)  d L_{(s -a)+}  \right] <\infty.
$$
Then the following  coupled linear  FBSDE 
admits a unique solution $\Theta (t):= (x(t), y(t), z(t)):$ 
  \begin{equation}\label{FBSDE-L}
\left\{\begin{aligned}
 x (t)=&\  x_0+ \int_0^t \( -y(s) + b_0(s) \)ds+  \int_0^t\(- y(s) + \delta_0 (s) \) d L_{(t-a)^+ } +\int_0^t \(-z (s)+\sigma _0(s)   \)  dB_{L_{(s-a)^+ }},\\
  y (t)=&\ x (T)+ \phi_0+ \int_t^T  \(  x (s) +g_0(s) \)ds +  \int_t^T \(  x (s)+  h _0(s) \)d L_{(s -a)+}    -\int_t^T z (s)dB_ {L_{(s -a)^+}} .\\
\end{aligned}
\right.
\end{equation}
  Moreover,  there is a constant $C=C(T)\geq 1$ depending only of $T$  so that 
   \begin{eqnarray} \label{5.1-1}
  && \EE \Big[ \sup\limits_{[0,T]}  \( x(t)^2+y(t)^2  \)    +\int_0^T  z(t)^2 d L_{(t-a)^+ }  \Big] \nonumber\\
  & \leq  &  C\EE \left[x_0^2+ \phi_0^2 +    \int_0^T  \( b_0(s)^2 +g_0 (s)^2 \) ds  + \int_0^T \( \de_0(s)^2 + h_0(s)^2 + \si_0^2   \)  d L_{(s -a)+}  
   \right].  
     \end{eqnarray}
 
\end{lemma}

\pf Suppose that  $(x(t),y (t),z (t))$ is a solution for \eqref{FBSDE-L}. Define  $\bar y (t) :=   y (t)-x(t)$. Note that  
 \[
 x (t) = x (T) +\int_t^T (y (s)  - b_0(s))ds + \int_t^T (y (s)  - \de_0(s))d L_{(s -a)+}   
  +   \int_t^T    \( z (s)- \si_0(s) \)d B_ {L_{(s -a)^+}}  . 
  \]
 \[ \bar y(t) = \phi_ 0 +\int_t^T\(-\bar y (s)+g_0(s)+b_0(s)\)ds +\int_t^T\(-\bar y (s)+h_0(s)+\de_0(s)\) d L_{(s -a)+} 
    -\int_t^T \( 2z_s-\si_0  \) dB_ {L_{(s -a)^+}} .    \]
 As a special case of  \cite[Theorem 4.5]{ZhangChen2023}, the following    BSDE:
  \begin{equation} \label{e:3.7}
  \bar y (t) = \phi_ 0 +\int_t^T\(-\bar y (s) +g_0(s)+b_0(s)\)ds +\int_t^T\(-\bar y  (s) +h_0(s)+\de_0(s)\) d L_{(s-a)+} 
    -\int_t^T \bar  z(s)  dB_ {L_{(s -a)^+}}    
   \end{equation} 
 has a unique solution  $(\bar y, \bar  z)$ that is $\{\sF'_t\}$-adapted.
 Define $z (t) :=\frac 1 2 \( \bar z (t)+\si_0(t)  \)$.
 Substituting $z (t) $ into the SDE for $x(t)$ in \eqref{FBSDE-L}    and noting that $y(t)= x(t)+\bar  y(t)$  yields:
 \begin{eqnarray}
 x (t) &= &  x_0+ \int_0^t \( -x (s) -\bar y (s) + b_0(s) \)ds+  \int_0^t\(-x (s) -\bar y (s)+ \delta_0 (s) \) dL_{(s-a)^+ }  \nonumber \\ 
 && +\int_0^t \(- \frac 1 2   \bar z (s)+ \frac 1 2   \si_0(s)  \)   dB_{L_{(s-a)^+ }}  . \label{e:3.8}
 \end{eqnarray} 
Clearly the above equation has a unique solution
\begin{eqnarray}
 x (t) &= & e^{-t-L_{(t-a)+}} x_0+ e^{-t-L_{(t-a)+}}  \int_0^t  e^{s+L_{(s-a)+}}  \(    b_0(s) -\bar y (s)   \)ds   \nonumber \\ 
 && +  e^{-t-L_{(t-a)+}}  \int_0^t  e^{s+L_{(s-a)+}} \(  \delta_0 (s)  -\bar y (s) \) dL_{(s-a)^+ } \nonumber \\
 && +\frac12 e^{-t-L_{(t-a)+}}  \int_0^t  e^{s+L_{(s-a)+}} \(\si_0(s) -   \bar z (s)    \si_0(t)  \)   dB_{L_{(s-a)^+ }}  . \label{e:3.8b}
 \end{eqnarray} 
 As $y (t):=\bar y (t)+ x (t) $, this proves that the  solution to \eqref{FBSDE-L} is unique.
Conversely, let $(\bar y (t), \bar z (t) )$ be the unique solution of the BSDE \eqref{e:3.7}.  Define $z (t)=\frac 1 2 \( \bar z (t)+\si_0(t)  \)$.
Let $x(t)$ be the unique solution of \eqref{e:3.8} and set $y (t)=\bar y (t)+ x (t) $. It is straightforward to verify that $(x(t), y(t), z(t))$
solves \eqref{FBSDE-L}. 

   We next proceed to establish  estimate \eqref{5.1-1}. 
It follows from \eqref{e:3.7} that $\tilde  y_s:= e^{-s-L_{(s-a)^+}}\bar y_s $ and  $\tilde  z_s:= e^{-s-L_{(s-a)^+}}\bar z_s $ satisfies
the BSDE
\begin{equation}\label{e:3.11} 
- d \tilde y_s: = e^{-s-L_{(s-a)^+}} \left (  (g_0(s)+b_0(s) ) ds + ( h_0(s)+\de_0(s))  d L_{(s -a)+}  \right)
    - \tilde  z_s dB_ {L_{(s -a)^+}}  
\end{equation} 
with $\tilde y_T= e^{-T-L_{(T-a)^+}} \phi_0$.  Define
$$
\xi:= e^{-T-L_{(T-a)^+}} \phi_0 + \int_0^T e^{-s-L_{(s-a)^+}}  (  (g_0(s)+b_0(s) ) ds  
+ \int_0^T  e^{-s-L_{(s-a)^+}}  ( h_0(s)+\de_0(s))  d L_{(s -a)+}   . 
$$
By Cauchy-Schwarz inequality,  there is a constant $C_1=C(T )>0 $ so that 
\begin{equation}\label{e:3.12} 
\E [ \xi^2] \leq C_1   \EE \left[\phi_0^2 +    \int_0^T  \(g_0 (s)^2  +  b_0(s)^2 \) ds  + \int_0^T \( h_0(s)^2 + \de_0(s)^2    \)  d L_{(s -a)+}   \right]. 
\end{equation}
By  \eqref{e:3.11},  for $t\in [0, T]$, 
\begin{eqnarray}\label{e:3.13a} 
\tilde y_t  &=&\EE \xi -  \int_0^t e^{-s-L_{(s-a)^+}}  (  (g_0(s)+b_0(s) ) ds  
- \int_0^t  e^{-s-L_{(s-a)^+}}  ( h_0(s)+\de_0(s))  d L_{(s -a)+}    \nonumber \\
&& + \int_0^t \tilde  z_s dB_ {L_{(s -a)^+}} . 
\end{eqnarray}
In particular 
$\int_0^T     \tilde  z_s dB_ {L_{(s -a)^+}}=  \xi - \EE \xi  $. Thus we have
$$
\E   \int_0^T  |\tilde  z_s|^2 d {L_{(s -a)^+}} = \EE [ (\xi -\EE \xi )^2] \leq \EE [ \xi^2].
$$
 It follows together with \eqref{e:3.12} that 
 \begin{eqnarray}\label{e:3.14} 
 \E   \int_0^T  |\bar   z_s|^2 d {L_{(s -a)^+}}  
\leq    C_1 \EE \Big[ \phi_0^2 +    \int_0^T  \(g_0 (s)^2  +  b_0(s)^2 \) ds  + \int_0^T \( h_0(s)^2 + \de_0(s)^2    \)  d L_{(s -a)+}  
 \Big]. 
\end{eqnarray}
Apply Ito's formula to \eqref{e:3.13a} and using the Burkholder-Davis-Gundy inequality, 
we have by \eqref{e:3.12} and  \eqref{e:3.14}  that for some $C_2=C_2(T, \kappa)\geq 1$, 
  \begin{equation} \label{e:3.15a} 
    \EE \Big[ \sup\limits_{[0,T]}  \bar y(t)^2       \Big]  
    \leq    C_2 \EE \left[\phi_0^2 +    \int_0^T  \(g_0 (s)^2  +  b_0(s)^2 \) ds  + \int_0^T \( h_0(s)^2 + \de_0(s)^2    \)  d L_{(s -a)+}   \right]. 
        \end{equation}
 Similarly, we have from \eqref{e:3.8b} together with \eqref{e:3.14}-\eqref{e:3.15a}  that for some $C_3=C_3(T, \kappa)\geq 1$,
 \begin{equation} \label{e:3.17}
   \EE \Big[  \sup\limits_{[0,T]}    x(t)^2 \Big]   \leq    
   C_2 \EE \left[x_0^2+ \phi_0^2 +    \int_0^T  \(g_0 (s)^2  +  b_0(s)^2 \) ds  + \int_0^T \( h_0(s)^2 + \de_0(s)^2  + \si_0^2   \)  d L_{(s -a)+}   \right]. 
     \end{equation}
     The desired estimate  \eqref{5.1-1} now follows from the fact that 
    $y(t)=x(t)+\bar y(t) $, $z(t)= \bar z(t)+\delta_0 (t)$ and \eqref{e:3.14}-\eqref{e:3.17}. 
This completes the proof of the lemma. 
\qed 

\medskip

    For $\al \in[0,1]$, define 
 \begin{eqnarray*} 
&&b^\al (t, x, y)\triangleq \al b(t, x, y)  -(1-\al) y,  \nonumber\\
&&\de^\al (t,  x, y, z)\triangleq \al  \de(t, x, y, z)  -(1-\al) y,  \nonumber\\
&&\si ^\al (t, x,y,z)\triangleq \al  \si (t, x,y,z )  -(1-\al) z,  \nonumber\\
&&h^\al (t, x,y,z)\triangleq \al  h(t, x,y,z ) -(1-\al) x,  \nonumber\\
&&g^\al (t, x,y )\triangleq \al g(t, x,y)  -(1-\al) x,  \nonumber\\
&&\phi^\al ( x )\triangleq \al \phi(x) + (1-\al) x.\end{eqnarray*}
The above functions satisfy all the conditions in Hypothesis \ref{HP0}
with  the monotonicity conditions \eqref{m1} and \eqref{m2}  with the constant  $c$ there 
being replaced by 
\begin{eqnarray} \label{cal} 
c_\al  := 
\al c+ (1-\al) \ge \min (c, 1) 
\end{eqnarray}
and with the same Lipschitz constant $L\geq 1$.

For notational simplicity, let  $\Gamma_s :=(x_s,y_s)$ and  $\Theta_ s :=(x_s,y_s,z_s)$, 
 and denote by  ${\rm FBSDE} (\al)$
 the following  FBSDE:
  \begin{equation}\label{L-al}
\left\{\begin{aligned}
 x (t)=&\  x_0+ \int_0^t \( b^\al(s,\Gamma_s) + b_0(s) \)ds+  \int_0^t\(\de^\al(s,\Theta_s)   + \delta_0 (s) \)dL_{(t-a)^+ }  \\
 &\ +\int_0^t \(  \si ^\al(s,\Theta_s) +\sigma _0(s)   \)  dB_{L_{(t-a)^+ }},\\
    y (t)=&\ \phi^\al( x_T)+\phi_0 + \int_t^T  \(  g^\al(s,\Gamma_s)   +g_0 \)ds
    +\int_t^T   \(  h ^\al(s,\Theta_s) +  h _0) \)d L_{(t -a)+}  \\ 
  &\    -\int_t^T z (t)dB_ {L_{(t -a)^+}} . 
\end{aligned}
\right.
\end{equation} 
We say ${\rm FBSDE} (\al)$ is solvable if \eqref{L-al} has a unique solution for any $\de_0,  \,\si_0 , \, g_0$ and  $h_0$ are $\{\sF'_t\}_{t\in [0, T]}$-progressively measurable 
 $L^2$-integrable processes on $\Omega \times [0, T]$ with respect to $\bP\times dt$ and for any $\phi_0\in L^2(\sF'_T)$.
 By Lemma \ref{L:3.4}, ${\rm FBSDE} (0)$ is solvable.

 \begin{lemma}\label{L:3.6}
Suppose  Hypothesis \ref{HP0} holds. 
  Suppose that $\rm FBSDE (\al_0)$ is solvable 
  for some $\alpha_0\in [0, 1)$
 and  that there is a constant $C_0=C_0(T, \kappa, L, c)\ge 1$,
  where $L$ is the  Lipschitz constant   for  $(b, \de, \si,  g, h, \phi)$ and   $c$ is  the monotonicity constant     in  
 \eqref{m1}-\eqref{m2},   so that any solution to 
$\rm FBSDE (\al_0)$ satisfies
\begin{eqnarray} \label{e:3.22}
  && \EE  \Big[ \sup\limits_{[0,T]}  \( x(t)^2+y(t)^2  \)    +\int_0^Tz(t)^2 d L_{(t-a)^+ }  \Big] \nonumber\\
  & \leq &   C_0 \bigg( x_0^2 +    \EE  \left[      \phi (0)^2 \right] +\EE [ \phi_0^2 ]+\EE  \int_0^T  \left(
   b (  t,0,0)^2  + b_0(t)^2  +  g  (  t,0,0)^2+ g_0 (t)^2 \right) dt  \\
           && \quad + \EE \int_0^T \left(    \de (  t,0,0,0  )^2 + \de_0(t)^2  + h ( t,0,0,0)^2  +h_0(t)^2 
   +   \si ( t,0,0,0  )^2 +  \si_0 (t)^2 \right)dL_{(t-a)^+}       \bigg) .   \nonumber 
          \end{eqnarray}
  Then there exists
 $ \eta_0 =\eta_0 (T, \kappa, L, c ) \in (0, 1]$ so that $\rm FBSDE (\al)$ is solvable for any $\al \in[\al_0, (\al_0+\eta_0)\wedge 1    ]$
 and the solution to $\rm FBSDE (\al)$    satisfies \eqref{e:3.22} with $C_0$ replaced by $2C_0$.
 \end{lemma}

\pf (i) (Existence.) For any  $\al \in [\al_0,  (\al_0+\eta_0)\wedge 1  ]$,
  where $\eta_0>0$ is a positive constant to be determined later,
    denote $ \eta  :=   \al-\al_0\leq \eta_0     $. 
  For any  $\{\sF'_t\}_{t\in [0, T]}$- progressively measurable
 $L^2$-integrable processes $b _0,  \delta_0, \si_0, h_0$ and $ g_0$ on $\Omega \times [0, T]$ 
  with respect to $\bP\times dt$ 
  and $\phi_0\in L^2(\cF'_T)$, we define 
$\Gamma^n_t =(x^n_t,y^n_t)$ and $\Theta^n_t=(x^n_t,y^n_t, z^n_t)$ recursively as follows.
Let  $\Gamma^0_t  := (0_0,0)$ and  $\Theta^0_t\ := (0,0,0)$ for any $t\in [0, T]$.
For $n\geq 0$, let $\Theta^{(n+1}_t:=(x^{n+1}_t,y^{n+1}_t, z^{n+1}_t)$ be the unique solution 
  in $\mathcal M_a [0,T] $
of the following coupled FBSDE
 \begin{equation}\label{deduc}
\left\{\begin{aligned}
 x_t^{n+1}=&\  x_0+ \int_0^t \( b^{\al_0}(s,\Gamma_s^{n+1} )+ b_0^n(s) \)ds+  \int_0^t\(\de^{\al_0}(s,\Theta_s ^{n+1}  )   + \delta_0^n (s)  \)
 d L_{(s-a)^+ }  \\
 &\ +\int_0^t \(  \si ^{\al_0}(s,\Theta_s^{n+1}  ) +\sigma^n _0(s)   \)  dB_{L_{(s-a)^+ }},\\
    y_t^{n+1}=&\ \phi^{\al_0}( x_T^{n+1})+\phi_0^n + \int_t^T  \(  g^{\al_0}    (s, \Gamma_s^{n+1}   )   +g_0 ^n(s)\)ds
      \\ 
  &\    +\int_t^T   \(  h ^{\al_0}(s,\Theta_s ^{n+1}  ) +  h _0^n(s)) \)d L_{(s-a)+}  -\int_t^T z  ^{n+1}     (t)dB_ {L_{(s -a)^+}} , 
\end{aligned}
\right.
\end{equation} 
 where 
 \begin{eqnarray*} 
&&b^n_0(t)   \triangleq  \eta \( y^n_t +b  (t,\Gamma_t^n) \)+   b_0(t),  \nonumber\\
&&\de^n_0(t)   \triangleq  \eta \( y^n_t +\de (t,\Theta_t^n) \)+   \de_0(t),  \nonumber\\
&&\si ^n_0(t)   \triangleq  \eta \( z^n_t +\si  (t,\Theta_t^n) \)+  \si_0(t) , \nonumber\\
&&h ^n_0(t)   \triangleq  \eta\( x^n_t +h  (t,\Theta_t^n) \)+  h_0(t),  \nonumber\\
&&g ^n_0(t)   \triangleq  \eta \( x^n_t +g  (t, \Gamma_t^n) \)+   g_0(t),\nonumber\\
&&\phi ^n_0  \triangleq  \eta \( -x^n_T +\phi (x^n_T) \)+   \phi_0 . 
\end{eqnarray*} 
 Note that under sub-linear growth condition on $\phi$, $\phi (x^n_T) \in L^2(\sF'_T)$ and so is $\phi^n_0$. 
Similarly,  by the Lipschitz assumption  in Hypothesis \ref{HP0}(i),  
  $b^n_0, \de^n_0,  \,\si^n_0 , \, h^n_0$ and  $g^n_0$ are $\{\sF'_t\}_{t\in [0, T]}$-progressively measurable
 $L^2$-integrable processes on $\Omega \times [0, T]$ with respect to $\bP\times dt$.

 \medskip

 Set $\Delta\Theta^n := (\Delta x^n, \Delta y^n, \Delta z^n):= \Theta^{n+1}-\Theta^n$.
 Note that 
  \begin{equation*}\label{D-al}
\left\{\begin{aligned}
 d\Delta x_t^{n }=&\   \bigg( \big( b^{\al_0}(t,\Gamma_t^{n+1} ) -  b^{\al_0}(t,\Gamma_t^{n} ) \big) 
 + \eta\big( \Delta y_t^{n-1}+ b (t,\Gamma_t^{n} ) -  b (t,\Gamma_t^{n-1} )\big)  \bigg)      dt\\
 &\  +   \bigg(   \big( \de^{\al_0}(t,\Theta_t^{n+1} ) -  \de^{\al_0}(t,\Theta_t^{n} ) \big) 
 + \eta\big( \Delta y_t^{n-1}  + \de (t,\Theta_t^{n} ) -  \de (t,\Theta_t^{n-1} )\big)         \bigg)dL_{(t-a)^+ }  \\
 &\   +  \bigg(   \big( \si^{\al_0}(t,\Theta_t^{n+1} ) -  \si^{\al_0}(t,\Theta_t^{n} ) \big) 
 + \eta\big(  \Delta z_t^{n-1}  +\si (t,\Theta_t^{n} ) -  \si (t,\Theta_t^{n-1} )\big)         \bigg)     dB_{L_{(t-a)^+ }},\\
 d \Delta  y_t^{n}=&\    - \bigg(   \big( g^{\al_0}(t, \Gamma_t^{n+1}  ) -  g^{\al_0}(t, \Gamma_t^{n}  ) \big) 
  + \eta\big( \Delta x_t^{n-1} +  g (t,\Gamma^{n} _t ) -  g (t,\Gamma_t^{n-1} )\big)         \bigg)           dt\\
 &\   -   \bigg( \big( h^{\al_0}(t,\Theta_t^{n+1} ) -  h^{\al_0}(t,\Theta_t^{n} ) \big) 
 + \eta\big( \Delta x_t^{n-1} +   h (t,\Theta_t^{n} ) -  h (t,\Theta_t^{n-1} )\big)  \bigg)          d L_{(t -a)+}  \\ 
  &\   +\Delta z_t^n  dB_ {L_{(t -a)^+}}  . 
\end{aligned}
\right.
\end{equation*}  
 By   it\^o's formula,   we have
 \begin{eqnarray}  \label{e:3.25a} 
&&  d \(\Delta x_t^{n } \Delta  y_t^{n}\) \nonumber \\
 &=&\bigg(  -   \big( g^{\al_0}(t, \Gamma_t^{n+1}  ) -  g^{\al_0}(t, \Gamma_t^{n} ) \big)  \Delta x_t^{n } +\big( b^{\al_0}(t,\Gamma_t^{n+1} ) -  b^{\al_0}(t,\Gamma_t^{n} ) \big)\Delta y_t^{n } \bigg)  dt  \nonumber \\
 &&+\bigg(  -\big( h^{\al_0}(t,\Theta_t^{n+1} ) -  h^{\al_0}(t,\Theta_t^{n} ) \big) \Delta x_t^{n }+  \big( \de^{\al_0}(t,\Theta_t^{n+1} ) -  \de^{\al_0}(t,\Theta_t^{n} ) \big)\Delta  y_t^{n}\bigg) d L_{(t-a)^+ }  \nonumber \\
 &&+ \Delta z_t^n \big( \si^{\al_0}(t,\Theta_t^{n+1} ) -  \si^{\al_0}(t,\Theta_t^{n} ) \big) d L_{(t-a)^+ } 
 +\eta  \bigg(   \bigg(-  \big( \Delta x_t^{n-1}     +    g (t,\Gamma_t^{n} ) -  g (t,\Gamma_t^{n-1})\big) \Delta x_t^{n }
 \nonumber \\
&&  +  \big( \Delta y_t^{n-1} +    b (t,\Gamma_t^{n} ) -  b (t,\Gamma_t^{n-1} )\big)  \Delta  y_t^{n}     \bigg) dt
 +    \bigg(  -    \big( \Delta x_t^{n-1} +   h (t,\Theta_t^{n} ) -  h (t,\Theta_t^{n-1} )\big) \Delta x_t^{n }
 \nonumber \\
 && +  \big( \Delta y_t^{n-1}  + \de (t,\Theta_t^{n} ) -  \de (t,\Theta_t^{n-1} )\big)  \Delta  y_t^{n}
 +  \big(  \Delta z_t^{n-1}  +\si (t,\Theta_t^{n} ) -  \si (t,\Theta_t^{n-1} )\big)   \Delta z_t^n     \bigg) d L_{(t-a)^+ } \bigg) 
 \nonumber \\
 &&  + \hbox{martingale}   .
 \end{eqnarray}
  It follows from  Hypothesis \ref{HP0} and  \eqref{cal}   that 
 \begin{eqnarray} \label{e:3.26}
  &&\EE \bigg[  \int_0^T \Big(  -   \big( g^{\al_0}(t, \Gamma_t^{n+1}  ) -  g^{\al_0}(t, \Gamma_t^{n}  ) \big)  \Delta x_t^{n } +\big( b^{\al_0}(t,\Gamma_t^{n+1} ) -  b^{\al_0}(t,\Gamma_t^{n} ) \big)\Delta y_t^{n } \Big)  dt \nonumber \\
 && \quad + \int_0^T \Big(  -\big( h^{\al_0}(t,\Theta_t^{n+1} ) -  h^{\al_0}(t,\Theta_t^{n} ) \big) \Delta x_t^{n }+  \big( \de^{\al_0}(t,\Theta_t^{n+1} ) -  \de^{\al_0}(t,\Theta_t^{n} ) \big)\Delta  y_t^{n}\Big) d L_{(t-a)^+ }  \nonumber \\
 && \quad + \int_0^T \Delta z_t^n \big( \si^{\al_0}(t,\Theta_t^{n+1} ) -  \si^{\al_0}(t,\Theta_t^{n} ) \big) d L_{(t-a)^+ } \bigg] \nonumber \\
  && \leq    -c_{\al_0 }  \EE \int_0^T   |\Delta\Gamma_t^n |^2   dt  -  c_{\al_0 }  \EE \int_0^T      |\Delta\Theta^n_t |^2
   dL_{(t-a)^+}.   
 \end{eqnarray} 
On the other hand, by  the Lipschitz condition and Cauchy-Schwarz inequality, there is a constant $C_1 >0$ depending only
on $L$ and $c$ so that 
  \begin{eqnarray*}
 && \eta\EE \bigg[ \int_0^T    \Big(-  \big( \Delta x_t^{n-1}     +    g (t,\Gamma_t^{n} ) -  g (t,\Gamma_t^{n-1})\big) \Delta x_t^{n } +  \big( \Delta y_t^{n-1} +    b (t,\Gamma_t^{n} ) -  b (t,\Gamma_t^{n-1} )\big)  \Delta  y_t^{n}    \Big) dt\\
 &&\qquad +  \int_0^T   \Big(  -    \big( \Delta x_t^{n-1} +   h (t,\Theta_t^{n} ) -  h (t,\Theta_t^{n-1} )\big) \Delta x_t^{n }
  +  \big( \Delta y_t^{n-1}  + \de (t,\Theta_t^{n} ) -  \de (t,\Theta_t^{n-1} )\big)  \Delta  y_t^{n}\\
 && \hskip 0.8truein +    \big(  \Delta z_t^{n-1}  +\si (t,\Theta_t^{n} ) -  \si (t,\Theta_t^{n-1} )\big)   \Delta z_t^n    \Big) d L_{(t-a)^+ }  \bigg]\nonumber\\
 &\leq &  \eta\EE\Big [\int_0^T\( |\Delta \Gamma _t^{n }   | ^2 + C_1  |\Delta\Gamma _t^{n-1 }  | ^2   \) dt \Big] 
 + \eta \EE\Big [\int_0^T \( |\Delta\Theta _t^{n }   | ^2 + C_1  |\Delta \Theta _t^{n -1}  | ^2   \)dL_{(t-a)^+}\Big] .
 \end{eqnarray*}
  Integrating \eqref{e:3.25a} in $t$ over $[0, T]$ and then taking expectation, we have
from the above two inequalities that   
\begin{eqnarray}\label{e:3.27} 
&& \EE \left[ \Delta x_T^{n } \Delta  y_T^{n}      -\Delta x_0^{n } \Delta  y_0^{n}  \right]   \\
&\leq &  
\EE \int_0^T \Big(
(\eta-c_{\al_0 }  ) |\Delta\Gamma_t^n |^2  + C_1 \eta |\Delta\Gamma _t^{n-1 }  | ^2  \Big) dt + \EE \int_0^T
 \((\eta-c_ {\al_0 }   )  |\Delta\Theta _t^{n }   | ^2 +C_1 \eta  |\Delta\Theta _t^{n-1 }  | ^2   \)  dL_{(t-a)^+} .  \nonumber 
 \end{eqnarray}
Note that $\Delta x_0^n =0$ and    since the function $\phi^{\al_0} (x)$ is non-decreasing in $x$, 
\begin{equation}\label{e:3.28}
 \Delta  x^n_T \Delta  x^n_T=  \Delta  x^n_T \(   \phi^{\al_0}(x_T^{n+1}) -\phi^{\al_0} (x_T^{n+1}) \) \ge 0.
\end{equation} 
We have by \eqref{e:3.27} and \eqref{e:3.28} that 
\begin{eqnarray*} 
(c_ {\al_0 }   -\eta )  \Big(  \EE \int_0^T  
 |\Delta\Gamma_t^n |^2  dt +   \EE \int_0^T  |\Delta\Theta _t^{n }   | ^2   dL_{(t-a)^+}  \Big) \nonumber\\ 
  \leq C_1 \eta  \Big( \EE \int_0^T  |\Delta\Gamma _t^{n-1 }  | ^2    +   \EE \int_0^T    |\Delta\Theta _t^{n-1 }  | ^2   dL_{(t-a)^+} 
  \Big).
  \end{eqnarray*}
Without loss of generality we may and do assume that the monotonicy constant $c$ in \eqref{m1}-\eqref{m2} to be no larger than 1.  Then $c_{\al_0} \ge c$ by \eqref{cal} and thus 
 \begin{eqnarray*} 
 (c  -\eta )  \Big(  \EE \int_0^T  
 |\Delta\Gamma_t^n |^2  dt +   \EE \int_0^T  |\Delta\Theta _t^{n }   | ^2   dL_{(t-a)^+}  \Big) \nonumber\\ 
  \leq C_1 \eta  \Big( \EE \int_0^T  |\Delta\Gamma _t^{n-1 }  | ^2    +   \EE \int_0^T    |\Delta\Theta _t^{n-1 }  | ^2   dL_{(t-a)^+} 
  \Big).
     \end{eqnarray*} 
 Now take 
 \begin{equation}\label{e:3.28}
 \eta_0= \frac c{  1+4 C_1 },
 \end{equation} 
  which is strictly less than $c$. 
 Since $\eta \in [0, \eta_0]$, $\frac{C_1 \eta}{c-\eta} \leq \frac{C_1 \eta_0}{c-\eta_0}=1/4$, we have from the above display that 
  \begin{eqnarray*} 
 &&\EE  \left[ \int_0^T 
  |\Delta\Gamma_t^n |^2 dt   +  \int_0^T      
      |\Delta\Theta _t^{n }   | ^2   dL_{(t-a)^+} \right] \nonumber\\ 
  &\leq&  \frac 14 \,  \EE \left[ \int_0^T      |\Delta\Gamma _t^{n-1 }  | ^2  dt  + \int_0^T    |\Delta\Theta _t^{n-1 }  | ^2 dL_{(t-a)^+} \right]  .  
    \end{eqnarray*}   
    By induction, we have for every $n\geq 1$, 
\begin{equation}\label{e:3.24} 
 \EE  \left[  \int_0^T  |\Delta\Gamma_t^n |^2  dt  + \int_0^T    |\Delta\Theta _t^{n }   | ^2  dL_{(t-a)^+} \right] 
 \le \frac  {1 }{4^n}  \EE  \left[  \int_0^T    |\Delta\Gamma_t^0 |^2 dt   +  \int_0^T    |\Delta\Theta _t^{0 }   | ^2  dL_{(t-a)^+} \right] 
  =: \frac  {C_2 }{4^n} .  
  \end{equation} 
  Thus $\{\Theta^n_t; n\geq 1\}$ is a Cauchy sequence in $(\mathcal M_a [0,T] , \| \cdot\|_{\mathcal M_a [0,T] })$
  and so it converges to some $\Theta_t =(x_t, y_t, z_t) \in \mathcal M_a [0,T] $.
Note that
$$
b^{\al_0 } (t,\Gamma^{n+1}_t)+ b_0^n(t )
=b^{\al_0 } (t,\Gamma^{n+1}_t)+ \eta \( y^n_t +b  (t,\Gamma_t^n) \)+   b_0(t) 
$$
converges in $L^2(d\mathbb P \times dt)$ as $n\to \infty$ to 
$$  b^{\al_0 } (t,\Gamma _t)+ \eta \( y _t +b  (t,\Gamma_t ) \)+   b_0(t) 
=   b^{\al_0+\eta } (t,\Gamma _t)+     b_0(t)\nonumber\\=b^\al(t,\Gamma_t)+b_0(t) .
$$
Other coefficients converge in a similar way. 
   Thus by \eqref{deduc},  the Burkholder-Davis-Gundy inequality and the Lipschitz condition in Hypothesis \ref{HP0}(i)  we conclude
\begin{equation}\label{e:3.25}
\lim_{n\to \infty} \EE \left[ \sup_{t\in [0, T]} \( | x^n_t - x_t|+ | y^n_t - y_t| \)^2 \right] =0;
\end{equation}
 Consequently, 
 $\Theta_t$ is a solution to $ \rm FBSDE(\al) $.

\medskip

(ii) (Uniqueness.) 
We next show that solutions to $ \rm FBSDE(\al) $ is unique for any $\al \in [\al_0, \al_0+\eta_0]$. 
Suppose that  $\Theta = (x, y, z)$ and $\tilde \Theta =(\tilde x, \tilde y, \tilde z)$ are two solutions of $ \rm FBSDE(\al) $,
where $\al \in (\al_0, \al_0+\eta_0]$. 
For simplicity, let  $\Gamma :=(x,  y)$ and $\tilde \Gamma :=(\tilde x, \tilde y)$. 
Define 
\begin{eqnarray*} 
&&\Delta \Theta  :=(\Delta x, \Delta y, \Delta z) :=  \Theta- \tilde \Theta,   \hskip 0.3truein 
\Delta\Gamma :=  \Gamma- \tilde \Gamma,  \nonumber\\
&&\Delta\varphi_t  :=  \varphi(t,\Theta_t)-\varphi(t, \tilde \Theta_t)  \quad \hbox{for } \varphi=b^\al \hbox{ and }  g^\alpha,\nonumber\\
&&\Delta \xi_t  :=  \xi(t,\Gamma_t)-\xi (t, \tilde \Gamma_t)   \quad \hbox{for } \xi=\delta^\al, \sigma^\al 
\hbox{ and }  h^\al,\nonumber\\
&&\Delta \phi^\al_T  :=   \phi^\al (x_T )-\phi^\al ( \tilde x_T).\end{eqnarray*}
Then for $t\in [0, T]$, 
  \begin{equation*}
\left\{\begin{aligned}
 \Delta x (t)=&\ \int_0^t \Delta b^\al_s ds+\int_0^t   \Delta \delta^\al _s d L_{(s-a)^+ } +\int_0^t   \Delta \sigma^\al_s  dB_{L_{(s-a)^+ }},\\
  \Delta y (t)=&\  \Delta \phi^\al+\int_t^T  \Delta g^\al_s ds  +\int_t^T   \Delta h^\al_s d L_{(s-a)^+}    -\int_t^T  \Delta z (t)dB_ {L_{(s -a)^+}}  ,\\
\Delta x_T =&\ 0 \quad \hbox{and} \quad  \Delta y_T =\Delta \phi^\al_T .
\end{aligned}
\right.
\end{equation*}
By  It\^o's formula,   
\begin{eqnarray*}
d (\Delta x_t\Delta y_t) & =&  \( -\Delta x_t\Delta g^\al_t+ \Delta y_t\Delta b^\al_t \)dt 
+ \( -\Delta x_t\Delta h^\al_t+ \Delta y_t\Delta \delta^\al_t  +  \Delta \si^\al_t\Delta z_t   \) d L_{(t-a)^+ } \nonumber\\
&&+  \(  \Delta x_t\Delta z_t+ \Delta y_t\Delta \si^\al_t \)d B_{L_{(t-a)^+}}.
\end{eqnarray*}
Integrating the above in $t$ from $0$ to $T$ and then taking expectation yields
\begin{eqnarray}\label{e:3.15} 
 \EE \left[ \Delta x_T \Delta y_T\right]   
&=& \EE \left[\int_0^T  \(-\Delta x_t\Delta g^\al_t+ \Delta y_t\Delta b^\al_t \)dt \right] \nonumber \\
&& + \EE \left[ \int_0^T \( -\Delta x_t\Delta h^\al_t+ \Delta y_t\Delta \delta^\al_t  +  \Delta \si^\al_t\Delta z_t   \) d L_{(t-a)^+ } \right].
\end{eqnarray}
Note that by Hypothesis \ref{HP0}(iii), $\phi (x)$ is non-decreasing so is $\phi^\alpha (x)= \alpha \phi (x) + (1-\alpha) x$.
Thus 
\begin{equation}\label{e:3.32}
\Delta x_T \Delta y_T = (x_T-\tilde x_T)(\phi^\al (x_T)-\phi^\al (\tilde x_T)\geq 0.
\end{equation}  
Hence  we have by \eqref{e:3.15} and \eqref{cal} that 
\begin{equation}\label{e:3.33}
 0\leq  -c_\al \int_0^T \(  |\Delta x_t|^2+  |\Delta y_t |^2   \)dt  -     c_\al   \EE\int_0^T  
  \(|\Delta x  (t)  |^2  + | \Delta y  (t)  |^2  + |\Delta z  (t) |^2\) dL_{t-a)^+}.
\end{equation}  
Since $\Gamma_t$ and $\tilde \Gamma_t$ are continuous processes,
it follows that $\PP$-a.s. that $\Gamma_t= \Gamma_t$ for every $t\in [0, T]$ and 
$\EE \int_0^T (z_t -\tilde z_t)^2 dL_{(t-a)^+}=0$.
This proves that the solution to $ \rm FBSDE(\al) $ is unique. 

\medskip

 (iii) (Estimate.)  We now show estimate \eqref{e:3.22} holds for the  solution $\Theta_t =(x_t, y_t, z_t) $ of $\rm FBSDE (\al)$
with $2C_0$ in place of $4C_0$.  For simplicity, denote the quantity 
\begin{eqnarray*}
 && x_0^2   +     \EE  \left[  \phi (0)^2 \right]
 + \EE [ \phi_0^2 ]+\EE  \int_0^T  \left(
   b (  t,0,0)^2  + b_0(t)^2  +  g  (  t,0,0)^2+ g_0 (t)^2 \right) dt  \\
           && \quad + \EE \int_0^T \left(    \de (  t,0,0,0  )^2 + \de_0(t)^2  + h ( t,0,0,0)^2  +h_0(t)^2 +
    \si ( t,0,0,0  )^2 +  \si_0 (t)^2 \right)dL_{(t-a)^+}  
\end{eqnarray*}
by $A$.  If $A=0$,  the solution $\Theta_t$ is identically zero so the estimate
holds trivially. So we assume $A>0$. 
 Since $\Theta_t^{n+1}$ is a solution to \eqref{deduc}, which is a $\rm FBSDE (\al)$ with 
$(b^n_0, \de^n_0,  \si ^n_0, h ^n_0,   g ^n_0,   \phi ^n_0)$ in place of  $(b_0, \de_0,  \si_0, h_0,   g_0,   \phi_0)$,
 by assumption \eqref{e:3.22}, 
\begin{eqnarray*}  
  && \left( \EE  \left[ \sup\limits_{[0,T]}  \( x^{n+1}(t)^2+y^{n+1}(t)^2  \)    +\int_0^Tz^{n+1}(t)^2 d L_{(t-a)^+ }  \right]  \right)^{1/2}\nonumber\\
  & \leq &  \sqrt{ C_0} \bigg( x_0^2   +  \EE  \left[    \phi (0)^2 \right]   
  + \EE [ (\phi^n_0)^2 ]+\EE  \int_0^T  \left(
   b (  t,0,0)^2  + b^n_0(t)^2  +  g  (  t,0,0)^2+ g^n_0 (t)^2 \right) dt  \\
           && \quad + \EE \int_0^T \left(    \de (  t,0,0,0  )^2 + \de^n_0(t)^2  + h ( t,0,0,0)^2  +h^n_0(t)^2  
   +   \si ( t,0,0,0  )^2 +  \si^n_0 (t)^2 \right)dL_{(t-a)^+}       \bigg)^{1/2}   \nonumber  \\
  &\leq &       2 \sqrt{ C_0A }   +    \eta \sqrt{C_0}  \left( \EE  \left[ (\phi (x^n_T) -\phi (0)-x^n_T)^2 \right] \right)^{1/2} \\
&& +    2 \eta \sqrt{C_0} (1+L)  \bigg(   \EE  \int_0^T  
  \Gamma^n  (t)^2  dt    + \EE \int_0^T   \Theta^n (t)^2 dL_{(t-a)^+}       \bigg)^{1/2}   \\
  &\leq &     2 \sqrt{ C_0A } +    \eta \sqrt{C_0}  \left( \EE  \left[ (\phi (x^n_T) -\phi (0)-x^n_T)^2 \right] \right)^{1/2} \\
  && +
   +  2 \eta \sqrt{C_0} (1+L)  \bigg(  (1+\kappa^{-1}) \EE  \int_0^T  
  \Gamma^n  (t)^2  dt    + \EE \int_0^T   z^n (t)^2 dL_{(t-a)^+}       \bigg)^{1/2} ,    
              \end{eqnarray*}
              where in the second to the last inequality we used triangular inequality and the Lipschitz assumption on the coefficients
              and we used \eqref{e:2.2} in the last inequality.  
Taking $n\to \infty$,  
 we have by  $\lim_{n\to \infty}  \| \Theta _t -\Theta_t^n   \|_{\cM_a[0,T]} =0$ and \eqref{e:3.25} that 
\begin{eqnarray*} 
  && \left( \EE  \left[ \sup\limits_{[0,T]}  \( x (t)^2+y (t)^2  \)    +\int_0^T z (t)^2 d L_{(t-a)^+ }  \right]  \right)^{1/2}\nonumber\\
      &\leq & 2 \sqrt{ C_0A } +    \eta \sqrt{C_0}  \left( \EE  \left[ (\phi (x _T) -\phi (0)-x_T)^2 \right] \right)^{1/2} \\
   && +  2 \eta \sqrt{C_0} (1+L)  \bigg(   \EE  \int_0^T  \left(  (1+\kappa^{-1})
  x(t)^2  +  y   (t)^2 \right) dt    + \EE \int_0^T   z  (t)^2 dL_{(t-a)^+}       \bigg)^{1/2} \\
    &\leq & \sqrt{ C_0A } +    \eta \sqrt{C_0}  (L+1) \left( \EE  \left[ ( x_T^2 \right] \right)^{1/2} \\
   && +  2 \eta \sqrt{C_0} (1+L)  \bigg(    (1+\kappa^{-1}) T  \, \EE \sup\limits_{[0,T]}  \( x (t)^2+y (t)^2  \) 
       + \EE \int_0^T   z  (t)^2 dL_{(t-a)^+}       \bigg)^{1/2} \\
       &\leq & \sqrt{ C_0A } +        3 \eta \sqrt{C_0} (1+L)  (  (1+\kappa^{-1})T + 1)  \bigg(     \EE \left[ \sup\limits_{[0,T]}  \( x (t)^2+y (t)^2  \)\right]
       + \EE \int_0^T   z  (t)^2 dL_{(t-a)^+}       \bigg)^{1/2} .
              \end{eqnarray*}
 By decreasing the value of $\eta_0$ previously set in \eqref{e:3.28} if needed, we now take   
$$
\eta_0 :=  \min \left\{ \frac1{3  (1+L)  ( (1+\kappa^{-1}) T+ 1)}, \ \frac{c}{1+4C_1} \right\}.
$$
We then have from the above display that for any $\alpha \in [\alpha_0, (\alpha_0+\eta_0)\wedge 1]$, 
$$
 \left( \EE  \left[ \sup\limits_{[0,T]}  \( x (t)^2+y (t)^2  \)    +\int_0^T z (t)^2 d L_{(t-a)^+ }  \right]  \right)^{1/2} 
 \leq 2 \sqrt{ C_0A }.
 $$
 This establishes the desired estimates. The proof of the lemma is now complete. 
  \qed

\medskip

\noindent {\bf Proof of Theorem \ref{exsi-uniqu}.}  By Lemma \ref{L:3.4}, ${\rm FBSDE} (0)$ has a unique solution
 and the estimate \eqref{e:3.22} holds for $\alpha_0=0$. 
Applying Lemma \ref{L:3.6} 
with $b_0(t)=\delta_0(t)= \sigma_0(t)= h_0(t)=g_0(t)=0$ 
$[T/\eta_0]+1$ times, we deduce immediately that  the corresponding ${\rm FBSDE} (1)$, which is
the  FBSDE \eqref{FBSDE}, 
 has a unique solution and it enjoys the desired estimate \eqref{e:3.7a}. 
\qed
 
\medskip

Analogous to Hypothesis \ref{HP0}, we can consider the following condition. 

\begin{hypothesis}\label{HP0b}
 \begin{enumerate}
 \item[\rm (i)] 
  $ \delta,  \, \si$ and  $h $ are uniformly Lipschitz continuous in $(x,y,z)$,    $b$ and $g$ are uniformly Lipschitz continuous in $(x,y)$, 
  and $\phi$ is uniformly Lipschitz continuous in $x$ with Lipschitz constant $L\geq 1$
    with
  $$
  \EE \left[ \delta^2 ({\bf 0}) +\si^2  ({\bf 0})+h^2 ({\bf 0})   + b^2 ({\bf 0})+g^2 ({\bf 0}) + \phi^2(0)\right] <\infty.  
	$$

 \item[\rm (ii)] 
There is some $c>0$  so that the following two monotonicity conditions hold
$\PP$-a.s.: 
for any $t>0$,  $x_1, x_2, y_1, y_2, z_1$ and $z_2$ in $\R$: 
\begin{eqnarray}\label{m1b}
  &&  \big(  b(t,x_1,y_1 )  -b(t,x_2,y_2  ) \big) (y_1-y_2  )  -
  \big(  g(t,x_1,y_1 )  -g(t,x_2,y_2 ) \big) (x_1-x_2  ) \nonumber\\
 &&  \quad \geq c \big( | x_1-x_2 |^2  + | y_1-y_2 |^2   \big) ,
\end{eqnarray}
and
\begin{eqnarray}\label{m2b}
&&    \big(  \si(t,x_1,y_1,z_1)  -\si(t,x_2,y_2,z_2) \big) (z_1-z_2  ) +  \big(  \de(t,x_1,y_1,z_1)  -\de(t,x_2,y_2,z_2) \big) (y_1-y_2  )   \nonumber\\
&& \qquad \  -\big(  h(t,x_1,y_1,z_1 )  -h(t,x_2,y_2,z_2 ) \big) (x_1-x_2  ) \nonumber\\ 
&&\quad \geq c \big( | x_1-x_2 |^2  + | y_1-y_2 |^2 + | z_1-z_2 |^2  \big).
\end{eqnarray}

\item[\rm (iii)] The function $\phi (x)$ is non-increasing in $x$.
    \end{enumerate}
\end{hypothesis}

\begin{remark} \label{R:3.10} \rm
By the same arguments with straightforward modifications in  the inequalities in \eqref{e:3.26}-\eqref{e:3.28} 
and in \eqref{e:3.32}-\eqref{e:3.33},
all  the results including Lemma \ref{L:3.6} and Theorem  \ref{exsi-uniqu}
of this section hold  under Hypothesis  \ref{HP0b}.  
\qed
\end{remark}

 \vskip 0.3truein
 
 {\small 
{\bf Shuaiqi Zhang}

\smallskip

     School of Mathematics, China  University of Mining and Technology,    
  Xuzhou,   Jiangsu, 221116, China.  
  
  \smallskip

        Email: \texttt{shuaiqiz@hotmail.com}

\bigskip
	
{\bf Zhen-Qing Chen}

\smallskip

    Department of Mathematics, University of Washington, Seattle,
WA 98195, USA.  

\smallskip

    Email:  \texttt{zqchen@uw.edu}

}

\end{document}